# Smallest Enclosing Sphere in 3D - Particle Swarm Optimization Approach

Netzer Moriya

siOnet Ltd. - Applied Modeling Research, Herzelia, 46445 Israel,
e-mail: netzer@si-O-net.com

November 7, 2023

## Abstract

We have employed Particle Swarm Optimization to address a stochastic variant of the Smallest Enclosing Sphere estimation problem. An efficient algorithm has been developed to ascertain the optimal center and radius of a sphere encompassing a cloud of points within a three-dimensional space. Our findings are benchmarked against simulated scenarios of the classical problem. Additionally, we elucidate several benefits of our proposed algorithm over Welzl's Algorithm.

# 1 Introduction

The problem of the *Smallest Enclosing Sphere* (SES), also known as the *Minimum Enclosing Ball* (MEB), is a foundational challenge in computational geometry. The task is to find the smallest d-dimensional sphere that encloses a finite set of points $P$ in $\mathbb{R}^d$. This requires identifying a center $c \in \mathbb{R}^d$ and a radius $r \in \mathbb{R}$, minimizing $r$, such that $\|p-c\| \leq r$ for all $p \in P$.

Several approaches, such as Welzl's recursive algorithm [1], Megiddo's deterministic algorithm for the planar case [2], the Random Sample Consensus algorithm [3] commonly used in computer vision, and various convex optimization techniques like the Interior-point method [4], have been proposed to approximate the smallest enclosing sphere (or ellipsoid) for a point set. These methods, along with heuristic algorithms for the SES problem, are reviewed in literature including [5, 6, 7, 8, 9].

However, when point clouds exhibit stochastic variations, it can be advantageous to relax the stringent requirement of encompassing *all* points. We thus introduce a modified SES problem, herein referred to as the mSES Problem, which seeks the smallest sphere that includes as many points as possible while minimizing the average distance of the excluded points to the sphere's surface.

This modified objective strikes a balance between enclosing a maximal number of points, minimizing the sphere's volume, and keeping the sphere proximal to the excluded points. For quantifying this proximity, we employ the Least Mean Squares (LMS) method, which we anticipate to be particularly useful in balancing the inclusion of points and volume minimization.

## 2 Modified Smallest Enclosing Sphere (mSES) Problem

### 2.1 Problem Statement

Given a set of *n* points in in three-dimensional space $\mathcal{P} = \{p_1, p_2, \ldots, p_n\} \subset \mathbb{R}^3$, the objective is to find the smallest enclosing sphere defined by its center $C$ and radius $r$ such that:

1. The sphere $S$ encloses as many points from $\mathcal{P}$ as possible,

2. The radius $r$ of the sphere $S$ is minimized,

3. For points outside the sphere $S$, the Least Mean Squares (LMS) error of their distances to the sphere's surface is minimized.

### 2.2 Formal Definition

Let:

- $\mathcal{P} \subset \mathbb{R}^3$, be the set of all given points (the *cloud*),
- $\mathcal{P}_{\text{inside}}$ be the set of points inside the sphere $S$,
- $\mathcal{P}_{\text{outside}}$ be the set of points outside the sphere $S$, i.e., $\mathcal{P}_{\text{outside}} = \mathcal{P} \setminus \mathcal{P}_{\text{inside}}$.

The LMS error for the points outside the sphere is defined as:

$$\text{LMS} = \frac{1}{|\mathcal{P}_{\text{outside}}|} \sum_{p \in \mathcal{P}_{\text{outside}}} (\|p - C\| - r)^2 \qquad (1)$$

Where $\|\cdot\|$ denotes the Euclidean norm.

The objective function to be minimized can be formulated as:

$$J(C, r) = -\lambda \times |\mathcal{P}_{\text{inside}}| + \alpha \times r + \beta \times \text{LMS} \qquad (2)$$

Where:

- $\alpha$ is a weight representing the importance of minimizing the sphere's radius,
- $\lambda$ is a weight representing the importance of maximizing the number of enclosed points,
- $\beta$ is a weight representing the importance of minimizing the LMS error for points outside the sphere.

The mSES problem can then be posed as:

$$\min_{C, r} J(C, r) \qquad (3)$$

Subject to:

$$\forall p \in \mathcal{P}_{\text{inside}}, \|p - C\| \leq r \qquad (4)$$

$$\forall p \in \mathcal{P}_{\text{outside}}, \|p - C\| > r \qquad (5)$$

The goal is to minimize the radius $r$ while ensuring that the majority of points are enclosed within the sphere. This problem is NP-hard [10], i.e., there is no known polynomial-time algorithm that can solve it exactly for arbitrary point sets, and involves nonlinear optimization process due to the norm constraints.

# 3 Solution Approach - Particle Swarm Optimization (PSO) for the mSES Problem

Particle Swarm Optimization (PSO) is a well-regarded evolutionary computation technique, introduced by Kennedy and Eberhart in 1995 [11]. PSO's inspiration comes from natural phenomena such as bird flocking and fish schooling, incorporating collective intelligence dynamics into the optimization process.

The algorithm operates with a swarm of particles, where each particle represents a potential solution to the optimization problem. The particles move within the solution space, influenced by their personal best positions and the global best known position, thus embodying both personal and social components of collective behavior.

Despite its successes in various domains, PSO's stochastic nature necessitates careful tuning of parameters to prevent premature convergence or entrapment in local optima. These challenges remain an active area of research, with ongoing developments in hybrid strategies and adaptive methods.

## 3.1 Justification

- The PSO algorithm is based on the mimicking of social behavior. Particles in the search space are influenced by their own experiences and the experiences of their neighbors. By constantly updating velocities and positions, particles seek out optimal or near-optimal solutions,

- The inertia weight *w* controls the momentum of the particles. A high inertia weight can speed up convergence but may overshoot optimal solutions. A lower inertia weight may lead to more refined exploration but slower convergence.

- The cognitive coefficient $c_1$ represents the importance of personal experiences, while the social coefficient $c_2$ represents the importance of the swarm's experiences,

- The objective function is multi-objective, balancing between enclosing as many points as possible, minimizing the radius of the sphere, and ensuring the Least Mean Squares error for points outside the sphere.

## 3.2 The Particle Swarm Optimization (PSO) - Algorithmic Description

Particle Swarm Optimization (PSO) exhibits its prowess through simple yet powerful iterative procedures. The following articulate the primary steps involved in the PSO algorithm.

### 3.2.1 Initialization

The initialization phase lays the foundation for the search process by generating a swarm of particles with random positions and velocities within the predefined bounds of the problem space. Formally, the $i$-th particle in the swarm is initialized as follows:

- Position, $\mathbf{x}_i \in \mathbb{R}^n$: A vector representing a potential solution in the $n$-dimensional problem space. Each component $x_{ij}$ is typically initialized randomly within the search space limits $[x_{\min}, x_{\max}]$,

- Velocity, $\mathbf{v}_i \in \mathbb{R}^n$: A vector that determines the particle's movement direction and step size. Each component $v_{ij}$ is usually initialized randomly within $[v_{\min}, v_{\max}]$,

- Personal best, $\mathbf{p}_i$: The best position encountered by the particle, initialized to its initial position $\mathbf{x}_i$.

The global best $\mathbf{g}$ is initialized as the position of the particle with the best fitness value in the initial swarm [11, 12].

### 3.2.2 Velocity and Position Update

At each iteration, the particles update their velocities and positions to search for the optimum. The velocity update incorporates three components: the inertia of the previous velocity, the cognitive component, and the social component. The velocity $\mathbf{v}_i$ and position $\mathbf{x}_i$ of particle $i$ are updated as:

$$\mathbf{v}_i(t+1) = w\mathbf{v}_i(t) + c_1 r_1 (\mathbf{p}_i - \mathbf{x}_i(t)) + c_2 r_2 (\mathbf{g} - \mathbf{x}_i(t)), \tag{6}$$

$$\mathbf{x}_i(t+1) = \mathbf{x}_i(t) + \mathbf{v}_i(t+1), \tag{7}$$

where:

- $w$ is the inertia weight controlling the impact of the previous velocity on the current velocity,

- $c_1$ and $c_2$ are positive constants called cognitive and social parameters, respectively,

- $r_1$ and $r_2$ are random numbers drawn from a uniform distribution in the interval $[0, 1]$.

The inertia weight $w$ plays a crucial role in balancing global and local search abilities and is often set to decrease over time to encourage convergence [13]. The cognitive term represents the particle's memory of its personal best position, while the social term represents the cooperative aspect of the swarm by guiding particles towards the global best position.

### 3.2.3 Stochastic Nature of PSO

The stochastic components in PSO are crucial for its search mechanism, enabling the algorithm to balance between exploration and exploitation in the solution space. This randomization is introduced in two primary ways: the initialization of particles and the stochastic variables within the velocity update equation.

During initialization, particles are scattered randomly across the solution space, which prevents the algorithm from being biased towards any particular region [12]. As the algorithm progresses, the velocity of each particle is updated by incorporating stochastic variables $r_1$ and $r_2$, which are uniformly distributed random numbers. These variables are multiplied by the cognitive and social scaling factors, $c_1$ and $c_2$, influencing the trajectory of the particle towards its personal best and the global best [14].

The stochastic nature of these updates allows particles to explore the search space effectively, potentially avoiding local optima and ensuring a diverse set of solutions. However, this randomness also means that the convergence of the algorithm is not deterministic, and different runs may produce different outcomes [15].

The successful application of PSO in various optimization problems depends on its stochastic dynamics, which can be tuned by adjusting parameters such as the inertia weight $w$, and the cognitive and social coefficients, $c_1$ and $c_2$. A careful balance between these parameters is necessary to achieve optimal performance [16].

### 3.2.4 Personal and Global Best Update

Following the velocity and position updates, each particle evaluates its new position's fitness. If the new fitness is better than the fitness of its personal best $\mathbf{p}_i$, the personal best is updated to the new position. Similarly, if a particle finds a position better than the current global best, the global best is updated to this new position. This allows the swarm to collectively move towards the most promising regions of the search space.

### 3.2.5 Mathematical Insights

The Particle Swarm Optimization (PSO) algorithm is deeply rooted in the principles of optimization theory and dynamical systems. This section explores the mathematical formalism underpinning PSO and the conditions that govern its convergence and stability.

### 3.2.6 Dynamical System Model of PSO

PSO can be conceptualized as a discrete-time dynamical system where the state of each particle is defined by its position and velocity. Let $\mathbf{X} \in \mathbb{R}^{n \times d}$ be the matrix representing the positions of all particles in an $n$-particle, $d$-dimensional swarm, and let $\mathbf{V} \in \mathbb{R}^{n \times d}$ be the corresponding velocity matrix. The PSO algorithm updates the state of the swarm according to the following iterative map:

$$\mathbf{V}(t+1) = w\mathbf{V}(t) + c_1\mathbf{R}_1 \odot (\mathbf{P} - \mathbf{X}(t)) + c_2\mathbf{R}_2 \odot (\mathbf{G} - \mathbf{X}(t)), \tag{8}$$

$$\mathbf{X}(t+1) = \mathbf{X}(t) + \mathbf{V}(t+1). \tag{9}$$

where $\mathbf{P}$ denotes the matrix of personal best positions, $\mathbf{G}$ denotes the global best position replicated across all particles, $\mathbf{R}_1$ and $\mathbf{R}_2$ are matrices of random values drawn from a uniform distribution, and $\odot$ represents the Hadamard product (element-wise multiplication).

### 3.2.7 Convergence Analysis

Convergence in PSO implies that all particles in the swarm stabilize at a particular point in the search space, which ideally corresponds to a global optimum of the objective function $f : \mathbb{R}^d \to \mathbb{R}$. From an optimization perspective, the PSO algorithm is a zeroth-order method, as it does not require gradient information of $f$.

A necessary condition for convergence is that the sequences $\{\mathbf{X}(t)\}$ and $\{\mathbf{V}(t)\}$ are bounded. This is achieved if the following condition holds:

$$0 < w < 1, \quad 0 < c_1 + c_2 < \frac{2}{(1-w)}, \tag{10}$$

which ensures that the coefficients of the dynamical system induce a contraction mapping for sufficiently large $t$ [17].

For a comprehensive review of the mathematical properties and convergence analysis of PSO, readers are directed to seminal works in the field [18, 19, 20].

### 3.2.8 Stability Conditions

The stability of the PSO algorithm can be examined using the concept of Lyapunov stability. A Lyapunov function $L : \mathbb{R}^d \to \mathbb{R}$ can be constructed such that $L(\mathbf{X}(t))$ decreases monotonically with $t$. One candidate for $L$ is the average squared distance to the global best position:

$$L(\mathbf{X}(t)) = \frac{1}{n} \sum_{i=1}^{n} \|\mathbf{x}_i(t) - \mathbf{g}\|^2. \tag{11}$$

If $L(\mathbf{X}(t + 1)) < L(\mathbf{X}(t))$ for all $t$, then the swarm is converging to the global best position $\mathbf{g}$, which is a stable fixed point of the dynamical system [18].

### 3.2.9 Advantages and Disadvantages

**Advantages**  Particle Swarm Optimization (PSO) has several distinctive advantages that contribute to its popularity:

- **Ease of Implementation**: PSO is conceptually simpler and easier to implement compared to other evolutionary algorithms. The update rules for particles are straightforward, requiring only basic linear algebra operations [21],

- **Fewer Parameters to Tune**: PSO generally requires tuning fewer parameters. The inertia weight, cognitive coefficient, and social coefficient are the primary parameters, in contrast to genetic algorithms, which require more complex parameter tuning [20],

- **Flexibility**: The algorithm can be easily adapted and extended to handle a wide variety of optimization problems, including discrete, continuous, and multi-objective problems [19],

- **Highly Parallelizable**: The evaluation of particle positions and the update of their velocities can be done in parallel, making PSO suitable for parallel and distributed computing environments [22],

- **Ability to Escape Local Optima**: PSO's stochastic nature and social sharing mechanism help particles to escape from local optima, which is a common issue in gradient-based optimization methods [23].

**Disadvantages**  Despite its benefits, PSO also exhibits certain limitations:

- **Risk of Premature Convergence**: Without proper parameter settings, PSO may converge prematurely to suboptimal solutions, especially in complex multimodal landscapes [24].

- **Dependency on Initial Conditions**: The performance of PSO can be highly dependent on the initial distribution of particles. A poor initial distribution may lead to slower convergence or failure to find the global optimum [18].

- **Sensitivity to Parameter Settings**: The choice of parameters can greatly affect the performance of PSO. Finding an appropriate balance between exploration and exploitation often requires empirical experimentation [25].

- **Scaling Issues**: PSO's performance can degrade as the dimensionality of the problem increases, necessitating modifications to the standard algorithm for high-dimensional problems [26].

- **Lack of Theoretical Guarantees**: The theoretical analysis of PSO is challenging, and thus, there are fewer guarantees regarding its convergence properties compared to certain other optimization methods [27].

## 3.3 Objective Function

The objective function $J(C, r)$ was previously defined as:

$$J(C, r) = -\lambda \times |\mathcal{P}_{\text{inside}}| + \alpha \times r + \beta \times \text{LMS} \tag{12}$$

Where:

$$\text{LMS} = \frac{1}{|\mathcal{P}_{\text{outside}}|} \sum_{p \in \mathcal{P}_{\text{outside}}} (\|p - C\| - r)^2 \tag{13}$$

## 3.4 PSO Iteration

For each iteration $t$:

1. For each particle $i$:

- Update its velocity using:

$$\mathbf{v}_i(t+1) = w \times \mathbf{v}_i(t) + c_1 \times \text{rand}() \times (\mathbf{pbest}_i - \mathbf{x}_i(t)) + c_2 \times \text{rand}() \times (\mathbf{gbest} - \mathbf{x}_i(t)) \quad (14)$$

- Update its position using:

$$\mathbf{x}_i(t+1) = \mathbf{x}_i(t) + \mathbf{v}_i(t+1) \quad (15)$$

- If $J(\mathbf{x}_i(t+1)) < J(\mathbf{pbest}_i)$:
  - Update $\mathbf{pbest}_i$ as $\mathbf{x}_i(t+1)$.
- If $J(\mathbf{x}_i(t+1)) < J(\mathbf{gbest})$:
  - Update $\mathbf{gbest}$ as $\mathbf{x}_i(t+1)$.

Discussing the results in a following section, we wish to compare the optimization method used (i.e. PSO) with the optimization process obtained by using the Welzl's algorithm It may thus be useful to describe the Welzl's algorithm and present its main characteristics.

## 3.5 Welzl's Algorithm

Welzl's algorithm solves the classic SES problem (see above) with expected linear time complexity, efficiently computing the minimal enclosing sphere for a given set of points in Euclidean space, making it one of the fastest algorithms for this purpose [1]. The method employs a recursive approach that builds the smallest enclosing sphere incrementally. Welzl's algorithm begins with a base case of an empty point set and recursively adds points. If a point is outside the current sphere, the algorithm finds a new sphere that includes this point while also enclosing the previous points. It elegantly utilizes the fact that the boundary of the minimal sphere is defined by either two, three, or, in the three-dimensional case, four points. The SES problem is inherently an optimization problem that can be mathematically formalized in the realm of computational geometry. Welzl's approach, with its expected linear time complexity, stands out for its efficiency and elegance [1].

The core insight behind Welzl's method is the characterization of the minimal enclosing sphere's boundary. In three-dimensional space, the minimal sphere can be uniquely determined by either two, three, or four points that lie on its surface. These points are known as support points, and they provide a simple yet powerful constraint that significantly reduces the complexity of the problem.

Mathematically, the algorithm is underpinned by the following properties:

- The minimal enclosing sphere of a set of points can be determined by considering only a subset of these points that lie on the sphere's boundary.

- If a set of points is enclosed by a sphere with radius $r$, then any subset of these points is also enclosed by a sphere with radius $r$ or smaller.

- For a recursive algorithm, it is sufficient to consider the case where the minimal enclosing sphere for a subset is defined by exactly three or four support points.

**Algorithmic Robustness**   The robustness of Welzl's algorithm is due to its deterministic outcome irrespective of the input's ordering. While the algorithm itself involves randomization, the expected linear time complexity holds for any permutation of the input points. This randomization is employed not to affect the solution's accuracy but to ensure efficient performance across different instances.

**Convergence and Stability**   Convergence in Welzl's algorithm is guaranteed by the fact that it operates within a finite search space with a clear stopping condition when all points are enclosed within the sphere, or when the sphere is defined by the maximum number of support points. The stability of the algorithm stems from its geometrically sound incremental construction approach, where each step maintains the invariant of the minimal enclosing sphere.

To enhance the original algorithm, several strategies can be implemented, such as:

- Heuristic pre-processing to identify a good initial subset of support points, potentially reducing the recursion depth.

- Post-processing to validate the minimal sphere against numerical inaccuracies, especially in floating-point arithmetic.

- Hybridization with other geometric algorithms to handle degenerate cases where multiple solutions are possible.

With these considerations, Welzl's algorithm not only provides a theoretically optimal solution to the SES problem but also exemplifies an algorithm that is robust, converges efficiently, and is stable against the variabilities inherent in computational geometry.

# 4 Results

We have applied the Particle Swarm Optimization algorithm method on two sets of experimental clouds in the light of the modified-*Smallest Enclosing Sphere* problem, as defined above. The first experimental data was constructed by a set of points around the origin (i.e., $C_a = (0, 0, 0)$), with the spread of points follows the surface area of a sphere of an arbitrary radius ($r_a$). A Gaussian stochastic component was added to each of the points. In the second experiment, we have added a similar sphere (i.e., $r_b = r_a$), to the first cloud, however shifted with respect to the first cloud, along the three Cartesian directions ($C_b = (r_a/2, r_a/2, r_a/2)$), with similar Gaussian stochastic component to each point of the second cloud. Figure 1 shows the two clouds used.

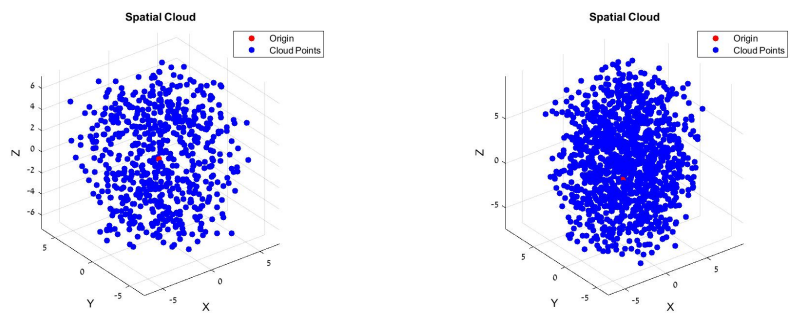

Figure 1: The first (left) and second (right) experimental clouds.

The optimal spheres according to the mSES problem modeled by the PSO algorithm are shown in figure 2. The first result (on the Left Hand Side (LHS)) refers to the optimal sphere for the first cloud (one sphere) while the second plot (Right Hand Side (RHS)) refers to the two-spheres cloud.

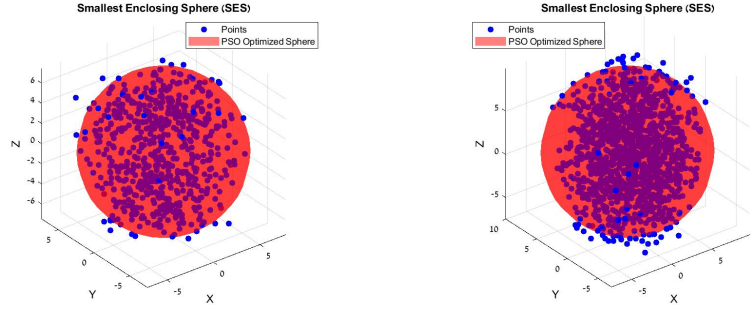

Figure 2: The first (left) and second (right) optimal spheres (Particle Swarm Optimization).

The optimal centers ($C_{1,2}^{PSO}$) and radii ($r_{1,2}^{PSO}$) were found to be:

$$C_1^{PSO} = (0.35, -0.05, -0.01) \quad r_1^{PSO} = 7.45$$
$$C_2^{PSO} = (1.08, 1.44, 1.23) \quad r_2^{PSO} = 8.66$$

for the first and second experiments respectively.

We have further compared the results to a classic calculation based on Welzl's algorithm [1] (see also the discussion above). Figure 3 shows the results obtained by the Welzl's algorithm where the first plot (on the LHS) refers to the optimal sphere related to the first cloud (one sphere) and the second (RHS) to the cloud build from the two spheres. Four points sets, are ensured to lie on each sphere optimized by the Welzl's algorithm process as expected.

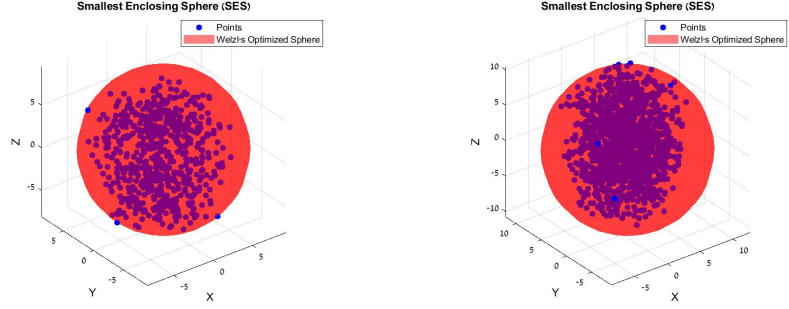

Figure 3: The first (left) and second (right) optimal spheres (Welzl's Algorithm).

The optimal centers ($C_{1,2}^{Wzl}$) and radii ($r_{1,2}^{Wzl}$) for this case were found to be:

$$C_1^{Wzl} = (0.44, -0.14, 0.71) \quad r_1^{Wzl} = 8.78$$
$$C_2^{Wzl} = (1.99, 1.58, -0.16) \quad r_2^{Wzl} = 10.57$$

The following can be commented:

1. In the first case (i.e., a single sphere), the optimization process based on the PSO-mSES algorithm yielded $C_1^{PSO}$ and $r_1^{PSO}$, with an optimal sphere that includes (in or on its surface) 601 point out of 634 (33 points left outside the sphere).

2. The corresponding result for the first case using the Welzl's algorithm showed $C_1^{Wzl}$ and $r_1^{Wzl}$. in this case, all points were included inside (or on) the sphere as expected.

3. In the second case (i.e., the two-sphere case, as described above), the optimization process based on the PSO-mSES algorithm yielded $C_2^{PSO}$ and $r_2^{PSO}$, with an optimal sphere that includes (in or on its surface) 1209 point out of 1268 (59 points left outside the sphere).

4. The corresponding result for the second case using the Welzl's algorithm showed $C_2^{Wzl}$ and $r_2^{Wzl}$. in this case, all points were included inside (or on) the sphere as expected.

5. As the Welzl's results strongly depend on the outmost points (still not to be considered as outliers), the PSO results are less sensitive to distant points (due to the LMS constraint).

# 5 Conclusions

The experimental results achieved by applying the Particle Swarm Optimization (PSO) algorithm to the modified Smallest Enclosing Sphere (mSES) problem present a compelling case for its use over Welzl's algorithm in scenarios where stochastic components influence the data points. We summarize the advantages as follows:

## 5.1 Advantages of PSO over Welzl's Algorithm

- **Robustness to Stochastic Data**: Unlike Welzl's algorithm, which is sensitive to the outermost points and may thus be influenced by outliers, the PSO algorithm shows robustness due to the inclusion of the Least Mean Squares (LMS) constraint. This allows PSO to minimize the influence of points that deviate significantly from the majority, resulting in a sphere that better represents the central tendency of the data.

- **Multi-Objective Optimization**: PSO inherently supports multi-objective optimization. The balance between enclosing as many points as possible, minimizing the sphere's radius, and minimizing the LMS error for the points outside the sphere presents a more nuanced approach to the mSES problem, which is beneficial in practical applications where complete enclosure is not strictly necessary.

- **Flexibility and Adaptability**: PSO's flexibility allows for easy adjustments to the objective function and constraints, making it adaptable to a wide variety of problems. This is particularly useful for mSES, where different problem instances may require different optimization strategies.

- **Parallelizability**: The PSO algorithm benefits from parallelizability, where the evaluation of each particle (potential solution) can be processed independently. This is advantageous over Welzl's recursive method, which is inherently sequential.

- **Efficiency in Higher Dimensions**: While Welzl's algorithm performs well in three dimensions, PSO is not limited by the dimensional complexity and can scale efficiently to problems involving higher-dimensional spaces.

## 5.2 Empirical Observations

- In the case of a single sphere, PSO produced a sphere with a slightly smaller radius than Welzl's algorithm, indicating a more efficient enclosure of the points.

- For the two-sphere cloud, PSO's performance was particularly notable. It achieved a balance between enclosing the majority of points and minimizing the sphere's radius, demonstrating its effectiveness in complex scenarios.

## 5.3 Conclusive Remarks

In conclusion, the application of the PSO algorithm to the mSES problem introduces a robust, adaptable, and efficient alternative to Welzl's algorithm, particularly when addressing the presence of stochastic elements in the data. PSO's ability to optimize multiple objectives concurrently proves invaluable in complex geometrical challenges, suggesting its suitability for a broader range of applications in computational geometry and related fields.